\documentclass[a4paper,11pt]{article}
\usepackage{amsfonts,amssymb,latexsym,amsmath}
\usepackage{multicol}
\newcommand{\la}{\lambda}
\newcommand{\gog}{{\mathfrak g}}
\newcommand{\nog}{{\mathfrak n}}
\newcommand{\mog}{{\mathfrak m}}
\newcommand{\bog}{{\mathfrak b}}
\newcommand{\ut}{{\mathfrak u}{\mathfrak t}}

\newcommand{\UT}{\mathrm{UT}}

\newcommand{\eps}{\varepsilon}

\newcommand{\Ad}{{\mathrm{Ad}}}
\newcommand{\Ab}{{\Bbb A}}
\newcommand{\Nb}{{\Bbb N}}
\newcommand{\BC}{{\cal B}}
\newcommand{\ZC}{{\cal Z}}
\newcommand{\AC}{{\cal A}}
\newcommand{\DC}{{\cal D}}
\newcommand{\LC}{{\cal L}}
\newcommand{\Dp}{\Delta_+}
\newcommand{\Mb}{{\Bbb M}}
\newcommand{\Zb}{{\Bbb Z}}
\newcommand{\Rb}{{\Bbb R}}
\renewcommand{\leq}{\leqslant}
\renewcommand{\geq}{\geqslant}

\begin{document}
\Large
\date{}
\title{Invariants of coadjoint representation of regular factors}
\author{A.N.Panov
\thanks{The paper is supported by RFBR grants 08-01-00151-a, 09-01-00058-a and by ADTP grant 3341}}
 \maketitle

\section*{\S 0. Introduction}

 Coadjoint orbits play an important role in the representation theory,
 symplectic geometry, mathematical physics. According to the orbit
 method of A.A.Kirillov ~\cite{K-Orb, K-62}, for nilpotent Lie groups there exists one to one correspondence
between coadjoint orbits  and irreducible representations in Hilbert
spaces.  This gives  possibility to solve problems of representation
theory and harmonic analysis in geometrical terms of the orbit
space. However, the problem of classification of all coadjoint
orbits for specific Lie groups (such as the group of unitriangular
matrices) is an open problem up today that is far from its solution.
In the origin paper~\cite{K-62} on the orbit method  the description
of algebra of invariants and classification of orbits of maximal
dimension was obtained.

 The main result of this paper consists in construction of
 generators  of the field of invariants for  the coadjoint representations of regular factors.
By regular factor we further call a Lie algebra that is a factor of
unitriangular Lie algebra with respect to  some regular ideal. The
paper consists of  three sections. In the first section we study the
diagramm method, introduced in  ~\cite{P1,P2}. In \S 2 we present
the notion of  extremal minor of characteristic matrix (see
definition 2.4). We shaw that its highest coefficient is invariant
with respect to the coadjoint representation (see theorem 2.5). The
method of proof is based on the reduction of quantum minors. We
state the conjecture 2.6 on structure of algebra of invariants
$K[\LC^*]^L$. In the last \S 3 we  prove that  the field of
invariants $K(\LC^*)^L$ is a field of rational functions on some
system of invariants (see theorem 3.20).

 Let $N=\UT(n,K)$ be the group of unitriangular matrices  of size  $n\times n$ with units on the diagonal and  with entries in the
 field $K$ of zero characteristic.
 The Lie algebra  $\nog=\ut(n,K)$ of this group consists of lower triangular matrices of size  $n\times n$ with zeros on the diagonal.
 One can define the natural representation of the group  $N$ in the conjugate space $\nog^*$ by the formula  $\Ad_g^*f(x)=f(\Ad_g^{-1} x)$,
 where  $f\in\nog^*$, $x\in\nog$ and $g\in N$.
 This representation is called   the coadjoint representation.
 We identify the symmetric algebra  $S(\nog)$ with the algebra of regular functions   $K[\nog^*]$ on the conjugate space  $\nog^*$.
Let us also identify $\nog^*$ with the subspace of upper triangular
matrices with zeros on the diagonal. The pairing  $ \nog$ and
$\nog^*$  is realized due to the Killing form  $(a,b)=
\mathrm{Tr}(ab)$, where  $a\in \nog$, $ b\in
 \nog^*$. After this identification the coadjoint action  may be realized by the formula  $\Ad_g^*b=P(\Ad_g b)$, where $P$ is the natural
 projection of the space  of  $n\times n$-matrices onto  $\nog^*$.

Recall that for any Lie algebra  $\gog$ the algebra  $K[\gog^*]$ is
a Poisson algebra with respect to the Poissson bracket such that
$\{x,y\}=[x,y]$ for any  $ x,y\in\gog$. In the case $k=\Rb$ the
symplectic leaves with respect to this Poisson bracket
 coincide with the orbits of coadjoint representation ~\cite{K-Orb}. Respectively,
the algebra  of Casimir elements in  $K[\gog^*]$ coincides with the
algebra of invariants  $K[\gog^*]^N$ of the coadjoint
representation.

The coadjoint orbits of the group  $N$ are closed with respect to
the Zariski topology in  $\nog^*$, since all orbits of a regular
action of an arbitrary algebraic unipotent group  in an affine
algebraic variety are closed  ~\cite[11.2.4]{Dix}.

To simplify  language we  shall give the following definition: a
root is an arbitrary pair $(i,j)$, where $i,j$ are positive integers
from 1 to $n$ and $i\ne j$. The permutation group  $S_n$ acts on the
set of roots by the formula  $w(i,j)=(w(i),w(j))$.

A root  $(i,j)$ is positive if  $ i > j$. Respectively, a root is
negative if  $i<j$. We denote the set  positive roots by $\Dp$.

For any root  $\eta=(i,j)$ we denote by  $-\eta$ the root  $(j,i)$.
We define the partial operation of addition  in the set of positive
roots: if $\eta=(i,j)\in \Dp$ and $\eta'=(j,m)\in \Dp$, then
$\eta+\eta'=(i,m)$.

Consider the standard basis  $\{y_{ij}:~ (i,j)\in \Dp\}$ in the
algebra
 $\nog$. We shall also use the notation $y_\xi$ for $y_{ij}$, where
 $\xi= (i,j)$.

An ideal  $\mog$ in the Lie algebra  $\nog$ is called regular, if it
is  generated by some subsystem of vectors of the standard basis.
Then $\mog= \mathrm{span}\{ y_\eta\vert~~ \eta\in M\}$, where $M$ is
a subset  of  $\Dp$, satisfying the following property: if in a sum
of two positive roots one of summands belongs to  $M$, then the sum
also belongs to $M$.

 Denote by  $\LC$ the Lie factor algebra  $\nog/\mog$ (the regular factor ) and by $L$
 the corresponding factor group of $N$ with respect to the normal subgroup
 $\exp(\mog)$.
Note that the conjugate space  $\LC^*$ is a subspace in $\nog^*$
which consists of  all $f\in\nog^*$ that annihilates  $\mog$. The
coadjoint $L$-orbit for  $f\in\LC^*$ coincides with its $N$-orbit.

\section*{\S 1. Diagram and permutation associated with Lie algebra   $\LC$ }

 In the paper  ~\cite{P1} we corresponded to any regular factor  $\LC$ the diagram $\DC_\LC$,
 constructed  applying some formal rule of arrangement of symbols
 in the table. By the diagram $\DC_\LC$ one can easily calculate the index of Lie algebra $\LC$.
Recall that the index of a Lie algebra is the minimal dimension of
centralizer  of a linear form on this Lie algebra. For algebraic Lie
algebras the index is equal to the transcendental degree of the
field of invariants of coadjoint  representation. For nilpotent Lie
algebras (for example, $\ut(n,K)$) the field of invariants of the
coadjoint representation is the pure transcendental extension of the
mail field of  degree being  the index~\cite{Dix}. Respectively, by
the diagram one can easily calculate the maximal dimension of
coadjoint orbits (see theorem  1.2). Earlier the diagram method was
used for classification of all coadjoint orbits of unitriangular
group of size $n\leq 7$ ~\cite{P3}, for description of special
families of coadjoint orbits for an arbitrary  $n$ (the subregular
orbits \cite{P3}; ~ the orbits, associated with
involutions~\cite{P4}).

Let us state the construction method  of the diagram $\DC_\LC$ and
formulate the main assertions of the papers ~\cite{P1, P2}. Consider
the  order  $\succ$ on the set  $\Dp$ such that
$$(n,1)\succ
(n-1,1)\succ\ldots\succ
 (2,1)\succ(n,2)\succ\ldots\succ(3,2)\succ\ldots\succ(n,n-1).$$

 By means of the ideal $\mog$ we construct the diagram that is a  $n\times
n$-matrix in that  all places  $(i,j)$, ~$i\leq j$, are not filled
and all other places (i.e. places of $\Dp$) are filled by the
symbols "$\otimes$"{}, ~"$\bullet$"{}, "$+$"\ and "$-$"\ according
to the following rules. The places  $(i,j)\in M$ are filled by the
symbol "$\bullet$"{}. We shall refer  the procedure of placing of
"$\bullet$" onto  the places  in $M$ as  the zero step in
construction of the diagram.

 We put  the symbol "$\otimes$"{} on the greatest (in the sense of order  $\succ$)  place  in $\Dp\setminus M$.
  Note that this symbol will take place in the first column if the set of pairs of the form $(i,1)$ in $\Dp\setminus M$ in not empty.
Suppose that we put  the symbol "$\otimes$"{} on the place $(k,t)$,~
$k>t$. Further, we put  the symbol "$-$" on all places  $(k,a)$, ~
$t< a <k$, and we put  the symbol "$+$"   on all places $(b,t)$,~
$1< b <k$. This procedure finishes the first step of construction of
diagram.

Further, we put  the symbol "$\otimes$"{} on the greatest (in the
sense of order $\succ$) empty place in  $\Dp$. As above, we put the
symbols "$+$"\ and "$-$"\ on  empty places, taking into account the
following:  we put  the symbols  "$+$"\ and "$-$"\ in pairs; if the
both places  $(k,a)$ and  $(a,t)$, where $ k>a>t$, are empty, we put
"$-$"\ on the first place and "$+$"\ on the second place; if one of
these places, $(k,a)$ or $(a,t)$, are  already filled, then we do
not fill the other place. After this procedure we finish the step
which we call the second step.

Continuing the procedure further we have got  the diagram. We denote
this diagram by  $\DC_\LC$. The number of last step is
equal to the number of the symbols "$\otimes$"{} in the diagram.\\
{\bf Example 1}. Let $n=7$, $\mog = Ky_{51}\oplus Ky_{61}\oplus
Ky_{71}\oplus Ky_{62}$. The corresponding diagram is as follows

\begin{center}
{\large $\DC_\LC$ =
\begin{tabular}
{|p{0.1cm}|p{0.1cm}|p{0.1cm}|p{0.1cm}|p{0.1cm}|p{0.1cm}|p{0.1cm}|}
\hline  &  &  &  & & &  \\
\hline $+$& & &  & & &  \\
\hline $+$&$+$ & & & & &  \\
\hline $\otimes$ & $-$ & $-$  &  & &  & \\
\hline $\bullet$ & $+$ & $+$&$\otimes$ & & &\\
\hline  $\bullet$ & $\otimes$ & $-$ &$+$  &$-$ & &\\
\hline   $\bullet$&  $\bullet$ & $\otimes$ &$\otimes$  &$-$  &$-$ & \\
\hline
\end{tabular}}
\end{center}

We construct this diagram in 5 steps, beginning  with the zero step:

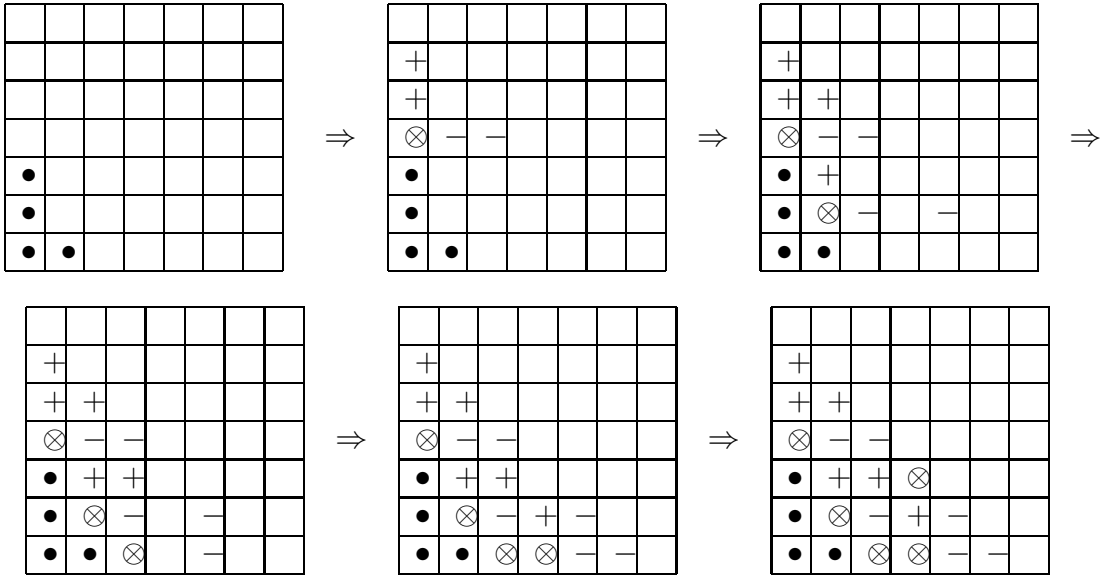
\begin{center}
{\large
\begin{tabular}{|p{0.1cm}|p{0.1cm}|p{0.1cm}|p{0.1cm}|p{0.1cm}|p{0.1cm}|p{0.1cm}|}
\hline  &  &  &  & & &  \\
\hline & & &  & & &  \\
\hline & & & & & &  \\
\hline  &  &   &  & &  & \\
\hline $\bullet$ &  & & & & &\\
\hline  $\bullet$ &  &  &  & & &\\
\hline   $\bullet$&  $\bullet$ &  &  &  & & \\
\hline
\end{tabular}
\quad $\Rightarrow$\quad
\begin{tabular}{|p{0.1cm}|p{0.1cm}|p{0.1cm}|p{0.1cm}|p{0.1cm}|p{0.1cm}|p{0.1cm}|}
\hline  &  &  &  & & &  \\
\hline $+$& & &  & & &  \\
\hline $+$& & & & & &  \\
\hline $\otimes$ & $-$ & $-$  &  & &  & \\
\hline $\bullet$ &  & & & & &\\
\hline  $\bullet$ &  &  &  & & &\\
\hline   $\bullet$&  $\bullet$ &  &  &  & & \\
\hline
\end{tabular}\quad $\Rightarrow$\quad
 \begin{tabular}{|p{0.1cm}|p{0.1cm}|p{0.1cm}|p{0.1cm}|p{0.1cm}|p{0.1cm}|p{0.1cm}|}
\hline  &  &  &  & & &  \\
\hline $+$& & &  & & &  \\
\hline $+$&$+$ & & & & &  \\
\hline $\otimes$ & $-$ & $-$  &  & &  & \\
\hline $\bullet$ & $+$ & & & & &\\
\hline  $\bullet$ & $\otimes$ & $-$ &  &$-$ & &\\
\hline   $\bullet$&  $\bullet$ &  &  &  & & \\
\hline
\end{tabular}\quad $\Rightarrow$}

\end{center}

\begin{center}
{\large
\begin{tabular}{|p{0.1cm}|p{0.1cm}|p{0.1cm}|p{0.1cm}|p{0.1cm}|p{0.1cm}|p{0.1cm}|}
\hline  &  &  &  & & &  \\
\hline $+$& & &  & & &  \\
\hline $+$&$+$ & & & & &  \\
\hline $\otimes$ & $-$ & $-$  &  & &  & \\
\hline $\bullet$ & $+$ & $+$& & & &\\
\hline  $\bullet$ & $\otimes$ & $-$ &  &$-$ & &\\
\hline   $\bullet$&  $\bullet$ & $\otimes$ &  &$-$  & & \\
\hline
\end{tabular}\quad $\Rightarrow$\quad
\begin{tabular}{|p{0.1cm}|p{0.1cm}|p{0.1cm}|p{0.1cm}|p{0.1cm}|p{0.1cm}|p{0.1cm}|}
\hline  &  &  &  & & &  \\
\hline $+$& & &  & & &  \\
\hline $+$&$+$ & & & & &  \\
\hline $\otimes$ & $-$ & $-$  &  & &  & \\
\hline $\bullet$ & $+$ & $+$& & & &\\
\hline  $\bullet$ & $\otimes$ & $-$ &$+$  &$-$ & &\\
\hline   $\bullet$&  $\bullet$ & $\otimes$ &$\otimes$  &$-$  &$-$ & \\
\hline
\end{tabular}\quad $\Rightarrow$\quad
 \begin{tabular}{|p{0.1cm}|p{0.1cm}|p{0.1cm}|p{0.1cm}|p{0.1cm}|p{0.1cm}|p{0.1cm}|}
\hline  &  &  &  & & &  \\
\hline $+$& & &  & & &  \\
\hline $+$&$+$ & & & & &  \\
\hline $\otimes$ & $-$ & $-$  &  & &  & \\
\hline $\bullet$ & $+$ & $+$&$\otimes$ & & &\\
\hline  $\bullet$ & $\otimes$ & $-$ &$+$  &$-$ & &\\
\hline   $\bullet$&  $\bullet$ & $\otimes$ &$\otimes$  &$-$  &$-$ & \\
\hline
\end{tabular}\quad\quad}
\end{center}

We denote  by $S=\{\xi_1\succ \xi_2\succ\ldots\succ \xi_s\}$ the set
of pairs  $(i,j)$, filled by  the symbol "$\otimes$"{} in the
diagram. For the diagram of example 1, the set  $S=\{\xi_1\succ
\xi_2\succ\xi_3\succ\xi_4\succ\xi_5\}$, where $ \xi_1=(4,1)$, $
\xi_2=(6,2)$, $ \xi_3=(7,3)$, $ \xi_4=(7,4)$, $ \xi_5=(5,4)$.

Denote by  $\Ab_m$  the Poisson algebra $K[p_1,\ldots,p_m;
q_1,\ldots,q_m]$ with the bracket $\{p_i,q_j\}= \delta_{ij}$.

Recall that a Poisson algebra  $\AC$ is a tensor product of two
Poisson algebras  $\BC_1\otimes\BC_2$, if $\AC$ is isomorphic to
$\BC_1\otimes\BC_2$ as commutative associative algebra  and   $\{\BC_1,\BC_2\}=0$.\\
 {\bf Theorem  1.1}~\cite{P1}. \emph{There exist  $z_1,\ldots,z_s\in K[\LC^*]^L$, where
$s=|S|$ such that  \\
1) any  $z_i=y_{\xi_i}Q_i+P_{>i}$, where $Q_i$ is some  product of
powers of    $z_1,\ldots, z_{i-1}$ and $P_{>i}$ is a polynomial in
variables  $\{y_\eta,~\eta\succ \xi_i\}$;
\\
2) denote by  $\ZC$ the set of denominators generated by
$z_1,\ldots,z_s$; the localization  $K[\LC^*]_\ZC$ of the algebra
$K[\LC^*]$ with respect to the set of denominators $\ZC$ is
isomorphic as a Poisson algebra to the tensor product
 $K[z_1^\pm,\ldots,z_s^\pm]\otimes \Ab_m$ for some  $m$.}
\\
Theorem 1.1 directly implies the following \\
{\bf Theorem  1.2}~\cite{P1}. \\ \emph{1) The field of invariants
$K(\LC^*)^L$ coincides with the field
$K(z_1,\ldots,z_s)$.\\
2) The maximal dimension of a coadjoint orbit in  $\LC^*$ equals to
the number of symbols
 "$+$"\ and "$-$"\  in the diagram  $\DC_\LC$.\\
 3) The index of Lie algebra  $\LC$ coincides with the number of symbols "$\otimes$"{} in the
 diagram $\DC_\LC$.}

To the Lie algebra $\LC$ we correspond the permutation,
 defined as follows.\\
{\bf Definition  1.3}~\cite{P2}. Denote by  $w = w_\LC$ the
permutation in
$S_n$ such that \\
1) ~ $w(1)=\max\{1\leq i\leq n|~ (i,1)\notin M\}$; \\
2) ~ $ w(t)=\max\{1\leq i\leq n|~ (i,t)\notin M,
~~i\notin\{w(1),\ldots,w(t-1)\} \}$  for all  $2\leq t\leq n$.

As usual, we denote by $l(w)$ the minimal number of  multipliers in
decompositions of  $w$ into products of simple reflections. The
number  $l(w)$ coincides with the number of inversions in the
rearrangement $(w(1),\ldots,w(n))$.
\\
{\bf Theorem 1.4}~\cite[Theorem 2.2]{P2}. \emph{The number $l(w)$ coincides with $\dim \LC$}.\\
{\bf Theorem  1.5}~\cite[Theorem  2.6]{P2}. \emph{We claim that $w =
r_{\xi_1}r_{\xi_2}\cdots r_{\xi_s}$}.

 Denote by  $\Dp^{(t)}$
the set  $\eta\in \Dp$ that have the form  $(b,t)$ for  $b>t$. Let
$S^{(t)} = \Dp^{(t)}\cap S$. Denote
$$S^{[t]} = S^{(1)}\sqcup\ldots\sqcup S^{(t)}.\eqno(1.1)$$

Denote by  $w^{(t)}$ (resp. $w^{[t]}$) the product of reflections
$r_\xi$ (arranged in the  decreasing order, in the sense of
$\succ$), where  $\xi\in S^{(t)}$ (resp. $S^{[t]}$).
 Easy to see that
$$w^{[t]}= w^{(1)}\ldots w^{(t)}.\eqno(1.2)$$
{\bf Theorem  1.6} ~\cite[Theorem  2.7]{P2}. \emph{Let  $\eta\in A^{(t)}$, then \\
1)~ the place  $\eta$ is filled in the diagram  $\DC_\LC$ by the
symbol  "$-$"\ iff
 $w^{[t-1]}(\eta)<0$;\\
2)~ the place  $\eta$ is filled in the diagram  $\DC_\LC$ by the
symbol
"$\bullet$"{} iff  $w^{[t]}(\eta)>0$;\\
3)~ место $\eta$ is filled in the diagram  $\DC_\LC$ by the symbol
"$+$"\ or "$\otimes$"{} iff $w^{[t-1]}(\eta)>0$ и $w^{[t]}(\eta)<
0.$}

For any  $\xi\in S^{(t)}$ we denote by  $w_\xi$ the product  of
reflections $r_{\xi'}$ (arranged in the  decreasing, in the sense of
$\succ$, order), where  $\xi'\in S$ и $\xi'\succeq \xi$.
 Easy to see that if  $\xi\in \Dp^{(t)}$, then
$$w_\xi= w^{(1)}\ldots w^{(t-1)} w_\xi^{(t)},\eqno(1.3)$$
where $w_\xi^{(t)}$ is the product  of reflections $r_{\xi'}$
(arranged in the  decreasing order, in the sense of $\succ$), where
$\xi'\in S^{(t)}$ and
$\xi'\succeq \xi$.\\
 {\bf Theorem  1.7}. \emph{Let $\xi =\xi_m \in S^{(t)}$ and $\eta\in \Dp$. Suppose that either
   a) $\xi \succ \eta$ or b) $\eta \in \Dp^{(t)}$. We claim that \\
1) if the place  $\eta$ was  empty after the  $m$th step or was
filled  by the  symbol   "$\bullet$"{}, then
$w_\xi(\eta)>0$;\\
2) if the place  $\eta$ was filled after the  $m$th step by any
symbol but not  "$\bullet$"{}, then
$w_\xi(\eta)<0$.}\\
{\bf Proof}.  If the place  $\eta$ is filled by the symbol
$\bullet$, then $w_\xi(\eta)>0$ (in the case  b) see
~\cite[Proposition 2.5]{P2}; in the case  a)  one can prove this
similarly. If the place $\eta$ is empty after the $m$th step, then
$w_\xi(\eta)>0$ (in the both cases  a) and b) see ~\cite[Proposition
2.3(1)]{P2}). If the place  $\eta$ was filled after the  $m$th step
by any symbol but not "$\bullet$"{}, then $w_\xi(\eta)<0$ (in the
case  a) see ~\cite[Proposition 2.3(2)]{P2}; in the case b) one can
prove similarly to ~\cite[Proposition 2.5(2)]{P2}).
 $\Box$\\
{\bf Corollary  1.8}. \emph{ For any  $\xi\in S$ the inequality  $w_\xi(\xi)<0$ holds.\\
{\bf Proof}. Let  $\xi = \xi_m$. Then $$w_\xi(\xi) =
w_{\xi_{m-1}}r_\xi(\xi) = -w_{\xi_{m-1}}(\xi)<0,$$ since the place
$\xi=\xi_m$ is empty after the  $(m-1)$th step.} $\Box$

\section*{\S 2. Reduction of quantum minors}

Recall some definitions of the theory of quantum matrices. Let  $q$
be a variable. The algebra of regular functions on quantum matrices
(briefly, the algebra of quantum matrices)  $K_q[Mat(n)]$ is
generated by the elements  $\{a_{ij}\}_{i,j=1}^n$ subject to the
system of relations  $ab=qba$, $cd=qdc$, $ac=qca$, $bd=qdb$,
$bc=cb$, $ad-da=(q-q^{-1})bc$, imposed on any  $2\times 2$-submatrix
$\left(\begin{array}{cc} a&b\\c&d\end{array}\right)$. The algebra
$K_q[Mat(n)]$ is a Noetherian ring without zero divisors and with
the Gelfand-Kirillov dimension equals $n^2$ (see for
instance~\cite{PW}).

 Denote by  $B:=B_-$ the group of lower triangular matrices and by  $\bog:=\bog_-$ its Lie algebra.
 One can construct the  algebra of regular functions  $K_q[B]$ on the quantum group
 $B$  factorizing   the algebra $K_q[Mat(n)]$  modulo the ideal   $<
 a_{ij}\vert~ i<j >$  and after  localizing with respect to the set of
 denominators,  generated by  $a_{11},\ldots, a_{nn}$.

The ideal  $\mog$ in Lie algebra Ли $\nog$ has its quantum analog --
the ideal  $Q_\mog$ in  $K_q[B]$, generated by  $a_{ij}$,~$(i,j)\in
M$. The factor algebra of  $K_q[B]$ modulo the ideal  $Q_\mog$ we
denote by  $K_q[L]$ and call the algebra of regular functions on the
quantum group  $L$.

For any systems of columns  $J=\{j_1<\ldots<j_m\}$ and rows
$I=\{i_1<\ldots<i_m\}$ the element
$$\Mb_I^J = \sum_{F=xI} (-q)^{l(x)} a_{f_1,j_1}\ldots
a_{f_m,j_m},\eqno(2.1)$$ where $F=(f_1,\ldots,f_m)$ and  $x\in S_m$,
is called the quantum minor.

 By definition, the quantum universal enveloping algebra  $U_q(\bog)$
 is generated by elements  $Y_1,\ldots, Y_{n-1}$ и $K_1,\ldots, K_n$ с
 subject to the relations  $$Y_iY_j^2 -(q+q^{-1})Y_iY_jY_i+Y_j^2Y_i = 0$$
 (the quantum Serre relations), where $ |
 j-i |=1$,  and
 $$K_iY_j=q^{-(\eps_i, \eps_j-\eps_{j+1})}Y_jK_i,$$ where  $1\leq i\leq n$,
 $1\leq j\leq n-1$.
By quantum algebra  $U_q(\nog)$ we mean the subalgebra  of
$U_q(\bog)$,  generated by $\{Y_i\}$. The ideal  $\mog$ corresponds
to the ideal  $\widetilde{Q}_\mog$ in $U_q(\nog)$, generated by all
$Y_{ij}$, where $i>j$ and  $(i,j)\in M$. We say that the quantum
group  $U_q(\LC)$ is the factor algebra of  $U_q(\nog)$ modulo the
ideal $\widetilde{Q}_\mog$.

Universal enveloping algebra  $U(\nog)$ is the factor algebra of
$U_q(\nog)$ modulo $q-1$. The symmetric algebra $S(\nog)$ coincides
with the graded algebra  $\mathrm{gr} U(\nog)$. Similarly, for the
Lie algebra Ли $\LC$ we have
$$S(\LC)=\mathrm{gr}\left(U_q(\LC)\bmod(q-1)\right).$$

It is obvious that the algebras   $ K_q[B]$ and $U_q(\bog)$ are not
isomorphic (since their factors modulo  $q-1$ are not isomorphic).
Denote by
 $ K'_q[B]$ and  $U'_q(\bog)$ the localizations of   $ K_q[B]$ and  $U_q(\bog)$ modulo $q-1$.

Well known that the algebras   $ K'_q[B]$ and  $U'_q(\bog)$ are
isomorphic. This isomorphism is called the isomorphism of Drinfeld.
 One can construct it directly: subalgebra of  $ K'_q[B]$, generated by the elements
 $$Y_{ij} = -\frac{a_{ij}a_{jj}^{-1}}{q-q^{-1}}, ~\mbox{where}~ i >
 j,\quad\mbox{and}\quad K_i = a_{ii}^{-1}, ~\mbox{where}~ 1\leq i\leq n, $$
is isomorphic to   $U'_q(\bog)$ and coincides with  $ K'_q[B]$ (see
for instance ~\cite{Ga}). Briefly, one can check this as follows: we
show that
  $\{Y_{i+1,i}, ~1\leq i\leq n-1\}$ satisfy quantum Serre relations; we extend the correspondence  $Y_i\to Y_{i+1,i}$
  to the homomorphism of  $U_q(\bog)$ into   $ K'_q[B]$;
this homomorphism is an isomorphism, since it induces the
isomorphism of corresponding graded algebras.

In what follows, we shall identify  $U_q(\bog)$ and $U_q(\nog)$ with
its images in  $ K'_q[B]$. Note that  the system of ordered (in the
sense of lexicographical order) monomials in $\{Y_{ij}\vert~ i>j\}$
forms the basis of  $U_q(\nog)$ as a free module over $K[q,q^{-1}]$.

 Consider the right action of  $n$-dimensional torus on the algebra of quantum matrices by the formular
 $a_{ij}\ldotp {\bf t} = t_ja_{ij}$, where ${\bf
t}=(t_1,\ldots,t_n)$.  We say that an element  $b$ of  $K_q[B]$ is
homogeneous of weight  $(k_1,\ldots, k_n)$, if
$$b\ldotp {\bf t}= t_1^{k_1}\ldots t_n^{k_n} b.$$
 By a homogeneous element
$b\in K_q[B] $ we construct the element  $$\tilde{b}=
\frac{(-1)^k}{(q-q^{-1})^k}\cdot b a_{11}^{-k_1}\cdots
a_{nn}^{-k_n}$$ in $U'_q(\nog)$, where $k=k_1+\ldots+ k_n$.

For the quantum minor  $\Mb_I^J$ of size  $m$ with system of columns
$I$ and rows $J$ we obtain
$$\widetilde{\Mb}_I^J =  \sum_{F=xI,~ f_\alpha >
j_\alpha \forall 1\leq \alpha\leq m} (-q)^{l(x)} q^{\phi(x)}
(q^{-1}-q)^{-\delta(F,J)}Y_{f_1,j_1}\ldots Y_{f_m,j_m},\eqno(2.2)$$
where $\phi(x)$ is some integer and $\delta(F,J) =
\mathrm{card}\{1\leq \alpha\leq m\vert~ f_\alpha = j_\alpha\}.$

We construct the formal matrix  $\Phi_\LC$  that all places in $M$
and also the places on  and upper the diagonal are  filled by zeros;
on the other places of the form  $(i,j)$,~ $i>j$ we put elements
$y_{ij}$ of the standard basis. For instance, for the Lie algebra
$\LC$ of the example 1 we obtain the diagram  $\DC_\LC$ and the
matrix $\Phi_\LC$:

\begin{tabular}{cc}{\large $\DC_\LC$ =
\begin{tabular}
{|p{0.1cm}|p{0.1cm}|p{0.1cm}|p{0.1cm}|p{0.1cm}|p{0.1cm}|p{0.1cm}|}
\hline  &  &  &  & & &  \\
\hline $+$& & &  & & &  \\
\hline $+$&$+$ & & & & &  \\
\hline $\otimes$ & $-$ & $-$  &  & &  & \\
\hline $\bullet$ & $+$ & $+$&$\otimes$ & & &\\
\hline  $\bullet$ & $\otimes$ & $-$ &$+$  &$-$ & &\\
\hline   $\bullet$&  $\bullet$ & $\otimes$ &$\otimes$  &$-$  &$-$ & \\
\hline
\end{tabular}}\quad,&

$ \quad\quad
 \Phi_\LC=\left(\begin{array}{ccccccc}
0&0&0&0&0&0&0\\
y_{21}&0&0&0&0&0&0\\
 y_{31}&y_{32}&0&0&0&0&0\\
  y_{41}&y_{42}&y_{43}&0&0&0&0\\
   0&y_{52}&y_{53}&y_{54}&0&0&0\\
    0&y_{62}&y_{63}&y_{64}&y_{65}&0&0\\
     0&0&y_{73}&y_{74}&y_{75}&y_{76}&0
\end{array}\right)
$
\end{tabular}

Let $\la$ be a variable. Any minor  $M_I^J$ of the matrix $\Phi_\LC$
is an element of $S(\LC)=K[\LC^*]$. Consider the characteristic
matrix  $\Phi_\LC - \la E$. The minor  $M_I^J(\la)$ of
characteristic matrix with  system of rows к $I$ and columns $J$ has
the form:
$$
M_I^J(\la)= \sum_{F=xI,~ f_\alpha > j_\alpha \forall 1\leq
\alpha\leq m} (-1)^{l(x)}  (-\la)^{\delta(F,J)}y_{f_1,j_1}\ldots
y_{f_m,j_m}. \eqno(2.3)
$$
Decomposing  (2.2) (resp. (2.3)) in powers of $(q-q^{-1})^{-1}$
(resp. $\la$), we obtain

$$ \widetilde{\Mb}_I^J = \sum_{\alpha=0}^n
\Mb_\alpha (q-q^{-1})^{-\alpha},\eqno(2.4)$$
$$M_I^J(\la) = \sum_{\alpha=0}^n M_\alpha\la^\alpha,\eqno(2.5)$$
where $\Mb_\alpha \in U_q(\nog)$,~ $M_\alpha \in U(\nog)$ for any
$\alpha$.

We shall say that the degree of quantum minor  ${\Mb}_I^J$ (more
precisely, the degree modulo $\mog$) is the greatest number $d$ such
that  $\Mb_\alpha\notin \widetilde{Q}_\mog$. The degree of
${M}_I^J(\la)$ is defined in the usual way.\\
{\bf Lemma 2.1}. \emph{Let  $\overline{Q}_\mog$ be the ideal of
$U_q(\nog)$, generated by  $\widetilde{Q}_\mog$ and $q-1$. We claim that\\
1)~ $M_\alpha =\mathrm{gr}\left( \Mb_\alpha\bmod
\overline{Q}_\mog\right)$,\\
2) the degrees of minors  ${\Mb}_I^J$ and $M_I^J(\la)$ coincide.}\\
{\bf Proof} of the statement 1) is obvious. Let us prove 2).
 According to PBW theorem, the monomials
$$Y_{f_1,j_1}^{k_1}\ldots
Y_{f_N,j_N}^{k_N},\eqno(2.6) $$ where $k_1,\ldots, k_N\in \Zb_+$ and
$(f_1,j_1)\succ\ldots\succ (f_N,j_N)$,  form the basis of
$U_q(\nog)$ over the field  $K$. The similar system of polynomials
in $\{y_{ij}\}$ form the basis of  $S(\nog)$ over $K$.

The quantum minor of (2.2) is presented as a linear combination of
basic monomials of the form  (2.6). Hence,  $\Mb_\alpha$ belongs to
 $\widetilde{Q}_\mog$ if and only if every  monomial, which is included as a summand in $\Mb_\alpha$, belongs to $\widetilde{Q}_\mog$.
In its turn, the monomial of (2.6) belongs to  $\widetilde{Q}_\mog$
if and only if it contains at least one element $Y_{f,j}$ of
$\widetilde{Q}_\mog$. Similar argumentation is true for any
coefficient  $M_\alpha$ from (2.5). Hence, $\Mb_\alpha$ belongs to
$\widetilde{Q}_\mog$ if and only if  $M_\alpha$ belongs to
$S(\nog)\mog$.~ $\Box$

Let  $d$  be the common degree of the minors  ${\Mb}_I^J$ and
$M_I^J(\la)$.
 Then
$$ M_{I}^J(\la) = M_d\la^d+M_{d-1}\la^{d-1}+\ldots+ M_0,\eqno(2.7) $$
where $M_d\ne 0$,  and the element  $\widetilde{\Mb}_I^J$, which is
taken modulo $\widetilde{Q}_\mog$, decomposes:
$$\widetilde{\Mb}_I^J=  (q-q^{-1})^{-d}\left( \Mb_d + (q-q^{-1})\Mb_{d-1} + \ldots +  (q-q^{-1})^d\Mb_0\right),\eqno(2.8)$$
where $\Mb_d\ne  0\bmod \widetilde{Q}_\mog.$

For a quantum minor  $\Mb=\Mb_I^J$ and a number  $1\leq i<n$ we
denote:
$$ \Mb^\downarrow =\left\{ \begin{array}{l}
\Mb_{(I\setminus i)\cup{i+1}}^J,~~ \mbox{if}~~ i\in I~~ \mbox{and}
~~ i+1\notin I,\\ 0, ~~ \mbox{in~~all~~other~~cases},
\end{array}\right.$$
$$ \Mb^\leftarrow =\left\{ \begin{array}{l}
\Mb_I^{(J\setminus i+1)\cup i},~~ \mbox{if}~~ i+1\in J ~~ \mbox{and}
~~ i\notin J,\\ 0, ~~ \mbox{in~~all~~other~~cases.}
\end{array}\right.$$
Similarly, we denote the minors
 $M^\leftarrow(\la)$ and $M^\downarrow(\la)$ for $M_{I}^J(\la)$.
The commutative relations of the algebra of quantum matrices imply
the following  \\
{\bf Lemma  2.2}. \emph{Let  $\Mb=\Mb_I^J$ and $a=a_{i+1,i}$, then
the following quantity takes place in the ring  $K_q[B]$:
$$
\Mb a -q^s a\Mb = -(q-q^{-1})\Mb^\leftarrow a_{i+1,i+1}
+(q-q^{-1})\Mb^\downarrow a_{ii}\eqno(2.9)$$ for some
$s\in\{0,1,-1\}$.}\\
{\bf Proof} directly follows from  ~\cite[Lemmas 4.1.5, 5.1.2]{PW}.
One can also prove the formula  (2.9), using the formula of $R$-matrix (see ~\cite[Chapter 7]{Ja}). $\Box$\\
 {\bf Definition  2.3}. A nonzero minor  $M_{I}^J(\la)$ of the characteristic matrix $\Phi_\LC-\la E$
 is extremal, if for any  $i$ the following inequalities hold
$\mathrm{deg} M^\leftarrow(\la) < \mathrm{deg} M(\la)$ and
$\mathrm{deg} M^\downarrow(\la) < \mathrm{deg} M(\la)$ (this
inequalities we also consider to be true  if  $M^\leftarrow(\la)=0 $
or $M^\downarrow(\la)=0$). Roughly speaking, a minor is extremal if
its degree decreases while
its rows are moving  down or its columns are moving to the left.\\
{\bf Remark 2.4}. Similarly, one can define an extremal minor of the
algebra of quantum matrices. From lemma 2.1, a minor $M_I^J(\la)$ is
extremal if and only if the corresponding quantum minor  $\Mb_I^J$
is extremal.

We denote by $P_{I}^J$ the highest coefficient $M_d$ in
decomposition (2.7) of the minor $M_{I}^J(\la)$ .
\\
{\bf Theorem  2.5}. \emph{The highest coefficient  $P_{I}^J$ of any
extremal minor is invariant with respect to the adjoint (resp.
coadjoint) representation of the group  $N$ in $S(\nog)$ (resp.  $K[\gog^*]$).}\\
{\bf Proof}. It is sufficient to prove that
$\mathrm{ad}_yP_{I}^J=0$, where $y= y_{i+1,i}$ and $1\leq i\leq
n-1$. Denote
$$Y = Y_{i+1,i} = -\frac{a_{i+1,i}a_{ii}^{-1}}{ (q-q^{-1})}.$$
The formula (2.9) implies
$$ \widetilde{\Mb} Y -q^s Y
\widetilde{\Mb} = \widetilde{\Mb^\leftarrow} -
\widetilde{\Mb^\downarrow}. \eqno(2.10)$$

Taking  (2.10) modulo $\widetilde{Q}_\mog$, we obtain

$$(q-q^{-1})^{-d}\left\{\left[\Mb_d+(q-q^{-1})\Mb_{d-1}+\ldots \right] Y -
Y \left[\Mb_d+(q-q^{-1})\Mb_{d-1}+\ldots \right]\right\} =$$
$$
(q-q^{-1})^{-d}\left\{\left[
\Mb^\leftarrow_d+(q-q^{-1})\Mb^\leftarrow_{d-1} + \ldots \right] -
\left[\Mb^\downarrow_d+(q-q^{-1})\Mb^\downarrow_{d-1}+\ldots
\right]\right\}.
$$
We cut down  $(q-q^{-1})^{-d}$. Since  $\Mb_I^J$ is an extremal
minor (see  remark  2.4), then  ~ $\Mb^\leftarrow_d
=\Mb^\downarrow_d =0 \bmod \widetilde{Q}_\mog$. Further, after
reduction modulo  $q-1$ we  obtain that  $\mathrm{ad}_y$ annihilate
 $\Mb_d\bmod(q-1)\in U(\nog)$.
Since  $$P_I^J=\mathrm{gr}\left(\Mb_d\bmod(q-1)\right),$$ then $\mathrm{ad}_yP_{I}^J=0$. $\Box$\\
 {\bf Conjecture 2.6}. \emph{The algebra of invariants  $K[\LC^*]^L$ of coadjoint representation of regular factor $\LC$ is generated
 by the
 highest coefficients of extremal minors.}

 The conjecture is true for the case  $\LC=\nog$ (i.e. $\mog=0$).
A corner minor  $M_i$ is a minor of the matrix  $\Phi$ that is lying
on the intersection of the first $i$ columns and last  $i$ rows.
 The algebra of invariants $K[\nog]^N $ is generated by the system of corner minors  $M_i$ (see ~\cite{K-62}), where
 $1\leq i\leq \left[ \frac{n}{2}\right]$, any of which is extremal.

\section*{\S 3. Field of invariants}

In this section we shall correspond to any  $\xi\in S$  an extremal
minor  $M_\xi(\la)$ (see definition  3.8). As it was shown in the
section  \S 2, its highest coefficient  $P_\xi$ is an invariant of
the coadjoint representation of $\LC$.

We shall prove that the field of invariants  $K(\LC^*)^L$
 coincides with the field of rational functions in  $P_\xi$, where $\xi\in S$.

Let  $\xi$ be some element of  $S$, say $\xi = (k,t)\in
S$, where $k>t$. As above,  $w_\xi$ is defined from  (1.3).\\
{\bf Lemma 3.1}. \emph{ Let  $i>t$. We claim that\\
1)  ~ if  $r_\eta(i)=i$ for any  $\eta\in S$ and $\eta\succeq \xi$,
then $w_\xi(i)=i$;
\\
2)~ if  $r_\eta(i)\ne i$ for some $\eta\in S$ and $\eta\succeq \xi$,
then
$w_\xi(i)\leq t$. In particular, $w_\xi(i)<i$; \\
3) ~$w_\xi(i)\leq i$.}\\
 {\bf Proof} of the statement  1) is obvious, the
statement  3) follows from  1) and 2). To prove 2) we consider the
sequence $i_0=i$, ~$i_1= w_\xi^{(t)}(i_0)$ and $i_\alpha=
w^{(t-\alpha+1)}(i_{\alpha-1})$, where
$2\leq\alpha\leq t$. \\
a)~ Under the assumption of statement 2), we shall show that there
exists a number  $1\leq\alpha\leq t$ such that
$i_\alpha=t-\alpha+1$.

Let $1\leq \alpha\leq t$  be the least number such that   the root
$\eta$, which is equal to  $(i,t-\alpha+1)$, contains in  $S$ and
$\eta\succeq \xi$. Then $i_{\alpha-1}=\ldots = i_0 = i$.

If set  $\{\theta\in S^{(t-\alpha+1)}|~ \theta\succ \eta\}$ is
empty, then  $i_\alpha =  r_\eta(i) = t-\alpha+1$. In the converse
case, let  $\theta=(j, t-\alpha+1) $ be the least, in sense of
$\succ$, element of  $S^{(t-\alpha+1)}$ such that  $\theta\succ
\eta$. Then  $j>i$ and $j= w^{(t-\alpha+1)}i_{\alpha-1} = r_\theta
r_\eta i_{\alpha-1} = i_\alpha$. The positive root  $\theta$ is the
sum of two positive roots  $\eta$  and  $\gamma=(j,i)$. Since
$r_\theta(\theta)<0$ (see corollary 1.8) and $r_\theta(\eta)>0$ (see
theorem 1.7(1)), then $r_\theta(\gamma)<0$. It follows from theorem
 1.7(2) that the place $\gamma$ will be filled by the symbol  "$-$" before
the  symbol "$\otimes$"{} appeared on the place  $\theta$. Hence,
there exists at least one  symbol "$\otimes$"{} in the row
$j=i_\alpha$ and in the columns with the numbers less or equal to  $
t-\alpha+1$. Let $\alpha'$ be the least number such that
$\alpha'>\alpha$ and the root  $\eta'$, which is equal to
$(i_\alpha,~ t-\alpha'+1)$, contains in $S$. Then
$i_{\alpha'-1}=\ldots = i_\alpha$.

If the set  $ \{\theta\in S^{(t-\alpha'+1)}|~ \theta\succ \eta'\}$
is empty, then $$i_{\alpha'} = r_{\eta'}(i_\alpha) = t-\alpha'+1.$$
In the converse case, the argumentation similar to above leads  to
existence of $\alpha''>\alpha'>\alpha$ such that  the root $\eta''$,
that is equal to $(i_{\alpha'},~ t-\alpha''+1)$, belongs to  $S$.
Continuing this process further, we obtain that at some $p$th step
the set
$$ \{\theta\in S^{(t-\alpha^{(p)}+1)}|~ \theta\succ \eta^{(p)}\}$$
is empty and, therefore,
$$i_{\alpha^{(p)}} = t-\alpha^{(p)} +1.$$ Note that for the constructed
$\alpha = \alpha^{(p)}$ we have inequality $i_\alpha=t-\alpha+1<t$
and $i_{\alpha-1}\geq
i_{\alpha-2}\geq\ldots \geq i$ at that.\\
 b) Let us finish the proof of statement  2) using induction on the number of elements in the set
 $\{\eta\in S\vert ~\eta\succeq\xi\}$. The proof is obvious if this set contains of one element.
 In the general case, for
$1\leq\alpha\leq t$ from  a) we obtain
$$w_\xi(i) = w^{(1)}\ldots w^{(t-\alpha)}w^{(t-\alpha+1)}\ldots
w^{(t-1)}w_\xi{(t)}(i) = w^{[t-\alpha]}i_\alpha =
w^{[t-\alpha]}(t-\alpha+1).$$ By induction assumption,
$w^{[t-\alpha]}(t-\alpha+1)\leq
t-\alpha+1 < t$. We conclude that $w_\xi(i)< t$. $\Box$\\
{\bf Corollary 3.2}. \emph{Let as above  $\xi=(k,t)$, ~$k>t$. We claim that \\
1)  ~ if  $r_\eta(k)=k$ for some  $\eta\in S$ and $\eta\succeq \xi$,
then $w_\xi(t)=k>t$;
\\
2)~ if  $r_\eta(k)\ne k$ for some   $\eta\in S$ and $\eta\succeq
\xi$, then
$w_\xi(t)< t$. } \\
{\bf Proof}. Since  $\xi\in S$, then  $\xi= \xi_m$ for some  $ 1\leq
m\leq s$ and $$w_\xi(t) = w_{\xi_m}(t) = w_{\xi_{m-1}}r_{\xi_m}(t) =
w_{\xi_{m-1}}(k).$$ Applying lemma 3.1 for $\xi=\xi_{m-1}$ and
$i=k$, we have got the proof of  statement of corollary.
$\Box$\\
{\bf Corollary 3.3}.\emph{ Let  $\xi$ be as in previous  corollary.
If $\xi$ is the least, in the sense of order
$\succ$, element of  $S^{(t)}$, then we claim that \\
1) ~ if  $r_\eta(k)=k$ for some $\eta\in S^{[t-1]}$ (i.e. there is
no symbol "$\otimes$"{}  in the $k$th row and columns with numbers
$< t$
of the diagram $D_\LC$), then $w(t)=k$;\\
\\
2)~ if  $r_\eta(k)\ne k$ for some  $\eta\in S^{[t-1]}$, then $w(t)<
t< k$.}
\\
 {\bf Proof}. For the least, in the sense of  $\succ$,  root $\xi\in S^{[t-1]}$
we have   $w_\xi=w^{[t]}$. The statement follows from corollary 3.2
applying  $w^{[t]}(j)=w(j)$ for any  $1\leq j\leqq t$.
 $\Box$\\
{\bf Corollary 3.4}. \emph{Any symbol  $\otimes $ of an arbitrary
$t$th
columns take place  either  $(w(t),t)$, or below it.} \\
{\bf Proof}. Let  $\xi$ be the least, in the sense of  $\succ$,
element of  $S^{(t)}$. If $r_\eta(k)=k$ for any  $\eta\in
S^{[t-1]}$, then $(w(t),t)=\xi$. If $r_\eta(k)\ne k$ for some
$\eta\in S^{[t-1]}$, then $w(t)<t$, the place   $(w(t),t)$ is upper
the diagonal and, therefore, it is upper all symbols "$\otimes$"{}.
~$\Box$

Let as above $\xi= ( k,t)$,~ $k>t$. Denote $$ h:= w_\xi(t).$$ We
shall give the definition of systems of columns  $J$ and rows $I$ of
the minor $M_\xi$ (see definition 3.8)  in each of the following
cases separately:~ 1) $h>t$ and
2) $h<t$. According to corollary  3.2, the case $h=t$ is impossible.\\
{\bf Case 1}. $h> t$. By corollary 3.2, $h=k>t$. In this case there
is no symbol "$\otimes$"{} in  the  $k$th row on the left side of
place
 $\xi=(k,t)$. We put
$$J:=J(\xi)=\{1\leq j\leq t:~w_\xi(j)\geq h\}, \quad\quad I:=I(\xi) =
wJ(\xi) .\eqno(3.1)$$ It is obvious that  $|I|=|J|$.\\
{\bf Case 2}. $h < t$. In this case there exists  at least one
symbol "$\otimes$"{} in the $k$th row on the left side of the place
$\xi=(k,t)$. The system  $J:=J(\xi)$ is defined as in (3.1). Denote
$$ I_*:=I_*(\xi)=\{t< i\leq n:~ i>t, ~w_\xi(i) < h\},$$
$$I:=I(\xi)=[h,t]\sqcup I_*.\eqno(3.2)$$
Here $[h,t]$ is a segment of positive integers (see definition
3.12(1)). Note that in the case  2 the equality  $|I|=|J|$ is not
obvious
beforehand and will be proved in lemma  3.7. \\
{\bf Remark 3.5}.\\
1) In both cases   $w_\xi(j) = w(j)$ for any  $1\leq j<t$. Hence,
$$ J= \{1\leq j<t:~~w(j)<h\}\sqcup \{t\}.\eqno(3.3)$$
2)~ Let  $\xi=\xi_m$. The element  ~$i\in I_*$ (see case 2) if and
only if  $w_\xi(i,t)<0$. By theorem  1.7,  the last  is equivalent
to the statement that the place  $(i,t)$ is filled after the $m$th
step by one of the symbols  "$\otimes$", "$+$"\ or "$-$" ( but not "$\bullet$"{}).\\
{\bf Lemma  3.6}.  \emph{Let as above  $\xi=\xi_m=(k,t)$. In both cases  1 and 2\\
1) there exists  $1\leq c\leq h$ such that  $J=[c,t]$;\\
2)  the rectangle  $[h,n]\times [1, c)$ in filled by the symbols
$\bullet$ in the diagram  $\DC_\LC$ (resp. by zeros in the matrix  $\Phi_\LC$);\\
3) there is no symbol "$\otimes$"{} in  the rectangle  $[1,h)\times J$ of the diagram  $\DC_\LC$.} \\
{\bf Proof}. Let  $w_\xi(j)<h$  for some  $1\leq j<t$. Since
$h=w_\xi(t)$, then  $h\notin\{w_\xi(1),\ldots, w_\xi(j-1)\}$. By
(3.3), $w_\xi(j) = w(j)< h$  and $h\notin\{w(1),\ldots, w(j-1)\}$.
The definition of  $w$ implies  $(h,j)\in M$ (i.e. the place $(h,j)$
is filled in the diagram by the symbol "$\bullet$"{}). Since $\mog$
is an ideal in  $\nog$, then  $ [h,n]\times [1,j]$ is contained in
$M$. This implies  1) and 2).

Since  $w(j)>h$ for some  $j\in [c,t)$, then, by corollary 3.4, all
symbols "$\otimes$"{}, which lie in columns  $[c,t)$, is contained
in rows
$[h,n]$. This proves 3). $\Box$\\
{\bf Lemma  3.7}.  \emph{Let  $\xi$ as above and  let  $h<t$ (i.e.
the case 2 takes place). Then
$|J|=|I|$. }\\
{\bf Proof}. Let  $a^{(t)}$ be the greatest number of  $[1,n]$ such
that  $(a^{(t)},t)$ do not belong to $M$. Then  $M\cap \Dp =
(a^{(t)},n]$. Easy to see that  $w_\xi(i)=i$ for any  $i\in
(a^{(t)},n]$. Therefore, $w_\xi[1, a^{(t)}] = [1, a^{(t)}]$.

Consider two decompositions of the segment  $[1, a^{(t)}]$ as a
union of disjoint sets:
$$[1, a^{(t)}] = [1,t) \sqcup [t, a^{(t)}],$$
$$[1, a^{(t)}] = [1,h) \sqcup [h, a^{(t)}].$$
Recall that $w_\xi(t)=h$. Let us show that

$$[1,h) = w_\xi[1,c)\sqcup w_\xi(I_*),\eqno(3.4)$$

$$ [h, a^{(t)}] = w_\xi(J)\sqcup w_\xi(I'_*),\eqno(3.5)$$
where $I'_*$ consists of all  $i$ such that the place $(i,t)$ was
not filled after the  $m$th step.

Since the  common number of elements of the left and right hand sets
of these formulas are equal, to prove (3.4) and (3.5) it is
sufficient to show that the right hand sets of these formulas are
contained in the corresponding right hand sets.

Really, if $i\in I_*$, then, by remark  3.5 and theorem 1.7 (case
b)), it follows  $w_\xi(i,t)< 0$. Hence,  $w_\xi(i)< w_\xi(t)=h$.
Therefore, $w_\xi(I_*)\subset [1,h)$.

If $i\in I_*'$, then, by theorem  1.7 (case b)), it follows
$w_\xi(i,t)>0$. Hence, $w_\xi(i)>w_\xi(t)=h$. That is
$w_\xi(I'_*)\subset [h, a^{(t)}]$.

The inclusion  $w_\xi(J)\subset [h, a^{(t)}]$ directly follows from
 definition of  $J$. Finally,
$w_\xi[1,c)\subset [1,h)$ follows from lemma 3.6(2). The formulas
(3.4) and (3.5) are proved.

By the equality (3.4), we have got   $ \vert w_\xi(I_*)\vert = h-c$.
Then $|I_*| = h-c$.  Since $ I = I_*\sqcup [h,t]$  and $ J = [c,t]$,
then
$$I=\vert I_*\vert + (t-h+1) = (h-c)+ (t-h+1)= t-c+1 = \vert
J\vert.~~\Box$$
 {\bf Definition 3.8}. In both cases  1) and 2) we
denote by $M_\xi(\la)$ the minor of characteristic matrix  $\Phi_\LC
-\la E $ with system of rows  $I=I(\xi)$ and system of columns
$J=J(\xi)$. We denote by $P_\xi$ its highest coefficient.

Our next goal is to prove that in both cases  $M_\xi(\la)$ is an
extremal minor; this will be proved in proposition 3.19. Case  1 is
more simple; reader interested in this case may turn directly to
proposition
 3.19. For the case  2 we need to prove the additional statements 3.9-3.18.

Let the case 2) takes place. As above  $\xi=\xi_m=(k,t)$, where
$k>t$, and $h=w_\xi(t)$, where $h<t$. As in the proof of lemma  3.7,
we define   $a^{(t)}:= \max\{i\vert ~ (i,t)\notin M\}$.
\\
{\bf Notations  3.9}.\\
1) ~ $E:= [h, a^{(t)}]$, ~~ $F_0:= I\setminus J = [c,h)$. \\
2) ~ $w_\xi^{[j]}:= w_\xi$, if $j\geq t$, and $w_\xi^{[j]}:=
w^{[j]}$, if $1\leq j<t$ (for definitions of  $w^{[j]}$ and  $w_\xi$ see (1.2) and (1.3)). \\
3)  Decompose  $E$ into the union of disjoint sets  $E=F\sqcup D$,
where
$$F:= \{a\in E|~ w_\xi^{[a-1]}(a)<h\},\eqno(3.6)$$
$$D:= \{a\in E|~ w_\xi^{[a-1]}(a)\geq h\}.\eqno(3.7)$$
4) ~$w_{\xi*}: = w^{(c)}\cdots w^{(t-1)}w_\xi^{(t)}$, $\quad\quad
w_{\xi*}^{[j]}:= \left\{\begin{array}{ll} w^{(c)}\cdots w^{(j)}&,
~\mbox{if}~
c\leq j< t,\\
w_{\xi*}&, ~\mbox{if}~ t<j\leq n.\end{array}\right.$\\

Note that  $w_\xi=w^{[c-1]}w_{\xi*} $. Since $w^{[c-1]}[1, h)\subset
[1,h) $  and $w^{[c-1]}(i)=i$ for any  $h\leq i\leq n$, then rewrite
(3.6) and (3.7) as follows

$$ F = \{j\in E|~ w_{\xi*}^{[j-1]}(j)\in F_0\},\eqno(3.8) $$
$$ D = \{j\in E|~ w_{\xi*}^{[j-1]}(j)\in J\}.\eqno(3.9) $$
{\bf Remark 3.10}. Note that  $h\in D$ (as there is no symbol
"$\otimes$" in the  $h$-row, otherwise  $w_\xi(t)<h$), ~$t\in F$ (as
$w_{\xi*}^{[t-1]}(t) = w_\xi(i_*)\in F_0$, where $i_*$ is an element
of  $I_*$ such that  $(i_*,t)$ is the greatest, in the sense of
$\succ$, element of  $S^{(t)}$) and $ I_*\subset
F$ (by definition of  $I_*$).\\
{\bf Lemma 3.11}. \emph{Let as above  $\xi=(k,t)$, where $k>t$. Let
$a\in D$,~ $b\in F$ and $a>b$. If $(b,p)\in
S$, where $p<t$,  then $(a,p)\in M$.}\\
{\bf Proof}. Let  $\xi_l = (b,p)\in S$. Let $q$ be the greatest
number such that  $q\leq m$ and $\xi_q\in\Dp^{[a-1]}$.
 Since
$b<a$, then $ l\leq q$.\\
1) Show that  $w_{\xi_i}(a)\geq h$ for any  $ 1\leq i\leq q$. We
shall prove by  induction on  $i$, starting with the greatest number
$q$. Since $a\in D$, then for  $i=q$ we obtain
$$w_{\xi_q}(a) = w_\xi^{[a-1]}(a)\geq h.$$ Suppose that
$w_{\xi_i}(a)\geq h$ was already proved for number $i$; let us prove
for number  $i-1$. Let $\xi_{i}=(a_1,b_1)$, ~$a_1>b_1$. If $a_1\ne
a$, then $r_{\xi_{i}}(a)= a$ and, therefore,
$$w_{\xi_{i}}(a) = w_{\xi_{i-1}}r_{\xi_{i}}(a) = w_{\xi_{i-1}}(a).$$
Using the induction assumption, we have got  $w_{\xi_{i-1}}(a)\geq
h$. If $a_1=a$, then
$$w_{\xi_{i}}(b_1) = w_{\xi_{i-1}}r_{\xi_{i}}(b_1) = w_{\xi_{i-1}}(a).$$
On the other hand, $w_{\xi_{i}}(\xi_{i})<0$ (see corollary 1.8).
Since $\xi_i =(a,b_1)$, then $$h\leq w_{\xi_{i}}(a)<
w_{\xi_{i}}(b_1)
= w_{\xi_{i-1}}(a),$$ this proves the statement 1).\\
2)~ Turn directly to the proof of the lemma. We may consider that
$p$ is a greatest number such that  $(b,p)\in S$, where  $p<t$.
Since  $b\in F$, then
 $$ w_{\xi_l}(b) = w_\xi^{[b-1]}(b) < h.$$
Introduce the notations  $\eta=(a,p)$,~ $\gamma = (a,b)$. We obtain
$\eta=\xi_l+\gamma$.\\
2a)~ Show that  $w_{\xi_{l}}(\gamma)>0$. Really, by  statement 1) of
the proof, it follows that  $w_{\xi_{l}}(a)\geq h$. On the other
hand, wee showed above that  $w_{\xi_{l}}(b) < h$. We obtain
$w_{\xi_{l}}(a)> w_{\xi_{l}}(b)$, hence  $w_{\xi_{l}}(\gamma) =
w_{\xi_{l}}(a,b) >0$.\\
2b)~ Since  $w_{\xi_{l}}(\gamma) >0$, then the place  $\gamma$ is
either empty after the  $l$th step, or filled by the symbol
"$\bullet$"{} (see theorem 1.7( case a)). If the last case takes
place, then  $(a,b)\in M$. Since  $\mog$ is an ideal of $\nog$, then
$(a,p)\in M$; this proves the statement of lemma.

Suppose that  $\gamma$ empty after the   $l$th step,  then the place
 $\gamma$ is also empty after the previous  $(l-1)$th step. By
 theorem 1.7 (case a),~ $w_{\xi_{l-1}}(\gamma)>0$. Since the place
$\xi_l$ is also empty after the  $(l-1)$th step, then
$w_{\xi_{l-1}}(\xi_l)>0$. Since  $\eta=\xi_l+\gamma$, then
$w_{\xi_{l-1}}(\eta)>0$. By  $a>b$, the place  $\eta$ is filled
after the  $(l-1)$th step. The theorem 1.7
implies that the place  $\eta$ maybe filled only by the symbol "$\bullet$"{} (i.e. $\eta\in M$). $\Box$\\
{\bf Notations and definitions 3.12}. \\
1)~ By a segment  $[a,b]$, ~$a,b\in \Nb$, of the positive integers
$\Nb$ we  call  the set  $\{ i\in \Nb|~ a\leq i\leq b\}$. \\
2) Let  $C\subset \Nb$. We shall say that a set  $A\subset C$ is a
segment in  $C$, if $A$ is an intersection of the positive integers
with  $C$. \\
 3)~ Introduce the relation $ < $ on the set of segments, according to which  $A<B$ if and only if
 $i<j$ for all  $i\in A$ and  $j\in B$.\\
4) If $A=[a,b]$ and $B=[c,d]$. Introduce the relation $A\lhd B$
(resp. $B\rhd A$) that means that  $a=c$ and $A\subseteq B$ (resp.
$b=d$ and $B\supseteq A)$.

Decompose  $D$ and $F$ into segments:
$$ D= D_1\sqcup D_2\sqcup \ldots \sqcup D_l, \quad \mbox{где}\quad D_1< D_2<\ldots < D_l,$$
$$ F= F_1\sqcup F_2\sqcup \ldots \sqcup F_l,  \quad \mbox{где}\quad F_1< F_2<\ldots < F_l$$
Note that   $$h\in  D_1<F_1< D_2< F_2 \ldots <D_l<F_l\rhd [k,
a^{(t)}].$$

 Let $i\in E$. Let  $\{i_\alpha\}$ be a sequence defined as follows. If  $i>t$, then, as in the proof of lemma  3.1, we put
  $i_0=i$, ~$i_1= w_\xi^{(t)}(i_0)$
and $i_\alpha= w^{(t-\alpha+1)}(i_{\alpha-1})$, where
$2\leq\alpha\leq t$. If $i\leq t$, then $i_\alpha$  is defined
similar, changing
$\xi$ by the least, with respect to  $\succ$, element of  $S^{[i-1]}$.\\
{\bf Definition 3.13}. For any  $i\in E$ we denote by  $i'$ the
greatest number in the sequence  $\{i_\alpha\}$, that is less than $
i$. By the proof of lemma 3.1 (see the end of a)), it  follows that
if $i'=i_\alpha$, then  $i_\alpha = i-\alpha+1 <i$ and
$i_{\alpha-1}\geq i_{\alpha-2}\geq\ldots\geq i$.

For any  $i\in I_*$ we construct a chain
$$i>i'>i''>\ldots> i^{(\mu(i))}=w_{\xi*}(i)\in F_0,\eqno(3.10)$$ in which
$i^{(p)}=(i^{(p-1)})'$ for any $0\leq p\leq \mu(i)$.\\
{\bf Lemma  3.14}. \\\emph{ 1) The chains  (3.10) for different
$i\in I_*$
do not intersect.\\
2) For any  $b\in F_0$ there exists  $i\in I_*$ such that
$b=i^{(\mu(i))}$.\\
3) If $a\in F\cap J$, then  $a=i^{(\nu)}$for some $i\in I_*$ and
$1\leq \nu\leq\mu(i)$.}\\
 {\bf Proof}. If  $a$ is a common element of chains (3.10) for
$i_1$ and $i_2$ from  $I_*$, then $$w^{(c)}\cdots w^{(a-1)}(a) =
w_{\xi*}(i_1)= w_{\xi*}(i_2).$$ Hence, $i_1=i_2$. This proves 1).

The statement 2) follows from  $w_{\xi*}(I_*)\subset F_0$ and
$|I_*|= h-c=|F_0|$ (see proof of lemma  3.7 and definition of $F_0$
from 3.9(1)).

Finally, let  $a\in F\cap J$ (i.e. $a\in F$ and $a\leq t$). Then $
w_{\xi*}^{[a-1]}(a)\in F_0$. There exists  $i\in I_*$ such that

$$w_{\xi*}^{[a-1]}(a) = w_{\xi*}(i).$$ We have  $a = w^{(a)}\cdots
w_\xi^{(t)}(i)$ and, therefore,  $a$ belongs to the chain  (3.10)
for $i$ (see the proof of lemma 3.7).
$\Box$\\
{\bf Definitions  3.15}.\\
1)
 For any $1\leq a\leq l$ we denote
$$D_{[a]} = D_1\sqcup\ldots \sqcup D_a,$$
$$F_{[a]} = F_1\sqcup\ldots
\sqcup F_a,$$
$$F'_a= \{i'|~ i\in F_a\},$$
$$F'_{[a]} = F'_1\sqcup\ldots
\sqcup F'_a.$$ 2) Note that  $k\in I_*\subset F_l$. Denote
$$F_{l1} = \{i\in F_l|~ i<k\},\quad\quad F_{l2} = \{i\in F_l|~ i\geq
k\},$$
$$ F'_{l1}= \{i'|~ i\in F_{l1}\},\quad\quad F'_{l2}= \{i'|~ i\in F_{l2}\} $$
3) Denote  $E_a=D_a\sqcup F_a$, for $1\leq a <l$,   and put
$E_l=D_l\sqcup F_{l1}$. Denote $$ E_{[a]}= E_1\sqcup\ldots\sqcup
E_a.$$ Note that $E=  E_{[l]}\sqcup F_{l2}$ and that $E_a$,
~$E_{[a]}$ and $F_0\sqcup E_{[a]}$ are segments of positive
integers, moreover  $E_{[a]}\lhd E$ and $F_0\sqcup E_{[a]}\lhd
[c,n]$ for any
$1\leq a\leq l$.\\
4)  For any  $1\leq a <l$ we denote by $J_a$ the subset
$F'_{[a]}\sqcup D_{[b(a)]}$, where $b(a)$ is the least number such
that  $F'_{[a]}\subseteq F_{[b(a)]}$ (in the case  $b(a)=0$ we put  $D_{[0]} = \varnothing$) .\\
5) By  $J_l$ we denote the subset  $ F'_{[l-1]}\sqcup F'_{l1}\sqcup
D_{[b(l)]}$, where  $b(l)$ is the least number such that  $F'_{[p-1]}\sqcup F'_{l1}\subset F_{[b(l)]}$.\\
{\bf Lemma  3.16}. \emph{
  We claim that  for any
$1\leq a\leq p$ \\
1) $F'_{[a]}\subset F_{[a-1]}$,\\
2) $F'_{[a]}$ is a segment in  $F$,\\
 3) $J_a$ is a segment in positive integers, moreover   $J_a\lhd F_0\sqcup E_{[b(a)]}$,\\
 4) ~ if  $i> E_{[a]}$ and $j\in J_a$, then $(i,j)\in M$,\\
 5) the minor of  matrix $\Phi_\LC$ with system of rows  $F_a$ and columns
 $F_a'$, where $1\leq a<l$, does not equal to zero,\\
 6) the minor of  matrix  $\Phi_\LC$ with system of rows $F_{l1}$ (resp. $F_{l2}$)   columns
 $F'_{l1}$ (resp. $F'_{l2}$)  does not equal to zero.}\\
{\bf Proof} will be carried out for   $l>a=1$. The general case can
be proved similarly, using induction by $a$.

For  $a=1$ we have  $F_{[1]} =F_1$  and  $J_1 = F'_{[1]} = F'_1$.
Let  $h_1$  be the least element of  $D_2$. Then $E_{[1]} =
[h,h_1)$.

Denote by  $g$  the greatest number such that $(i,g)\in S$ for some
$i\in F_1$. Then  $F'_1\subset [c,g]$. By lemma  3.11, ~$(h_1,g)\in
M$. Since  $\mog$ is an ideal in  $\nog$, then $[h_1,n]\times [1,
g]\subset M$, this proves  4).

Let us show that
$$F_1'=[c,g]\lhd  F_0.\eqno(3.11)$$ This implies the statements 1), 2) and 3) for the case  $a=1$.
\\
a)~ Show that  $g<h$ (recall that  $h$ is greatest element of
$D_1$). In the converse case,  $g\in D_1$ or $g> D_1$.

Suppose that  $g\in D_1$. Let  $(i,g)$ be the greatest, in the sense
of  $\succ$, element of  $S^{(g)}$ such that  $i\in F_1$. Then
$w^{(g)}(i) = g\in D_1$ and, therefore, $w^{[i-1]}(i)=
w^{[g-1]}w^{(g)}(i) = w^{[g-1]}(g) \geq h,$ this contradicts to
$i\in F$.

Suppose that  $g> D_1$. Then  $w^{(g+1)}\cdots
w^{(t-1)}w^{(t)}_\xi(t) \geq h_1$.  Since  $[h_1,n]\times [1,
g]\subset M$, then $w^{(c)}\cdots w^{(g)}(p)=p$ for any  $p\geq h_1$
and, therefore,
$$h= w_\xi(t)=w^{(c)}\cdots w^{(g)}w^{(g+1)}\cdots w^{(t)}_\xi(t)\geq
h_1.$$  On the other hand,  by definition  $h< h_1$. A contradiction. We obtain  $g< h$.\\
b)~ Conclude the proof of  (3.11). On one hand,  $w_{\xi*}(I_*) =
F_0=[c,h)$. On the other hand,
 $$w^{(g+1)}\cdots
w_\xi^{(t)}(I_*)\subset \left(F_0\setminus [c,g]\right)\sqcup
F_1\sqcup F_{\geq 2},$$ where $F_{\geq 2}$ is a union of all $F_p$,
~ $p\geq 2$, and $$w^{(c)}\cdots w^{(g)}(i)= i$$ for all  $i\in F_0$
 (see lemma 3.6(3)) and all  $i\in F_{\geq 2}$( see a)). Hence,
$$[c,g] =
w^{(c)}\cdots w^{(g)}(F_1)= F_1',$$ this proves  (3.11).\\
c) Let us prove the statements 5). Let  $f_1=|F_1|$. By the formula
(3.11), $g=c+f_1-1$. Let us show that there exist the systems
$\{i_{(p)}|~ 0\leq p\leq f_1-1\}$ and $\{j_{(p)}|~ 0\leq p\leq
f_1-1\}$ such that  $ (i_{(p)}, j_{(p)})\in S$ for all  $0\leq p\leq
f_1-1$ and
$$ F_1 = \{ i_{(0)},\ldots, i_{(f_1-1)}\}, \quad\quad  F'_1 = [c,g] = \{ j_{(0)},\ldots,
j_{(f_1-1)}\}.$$ Note that for an arbitrary $i\in F_1$ and  $j\in F_1'=[c,g]$ we have:\\
I) ~ $w^{(j)}(i) = i$, if $i\in F_0$ (see lemma  3.6(3));\\
II) ~ $w^{(j)}(i) = j\in F_0$, if $i\in F_1$ and  $(i,j)$ is the
greatest, in the sense of  $\succ$, element of  $S^{(j)}$;\\
III) ~ $w^{(j)}(i) \in F_1$, if $i\in F_1$  and $i$ does not satisfy
 II).

 Put $F_1=F_{10}$. For any  $0\leq p\leq f_1-1$ the following decomposition takes place
$$w^{(c+f_1-p)}\cdots w^{(c+f_1-1)}(F_1) = \{c-f_1-p,\ldots ,
c+f_1-1\}\sqcup F_{1p},\eqno(3.12) $$ where $F_{1p}$  is some subset
of  $F_1$ and $|F_{1p}|=f_1-p$. Since
$$F_1=F_{10}\supset F_{11}\supset\ldots\supset F_{1f_1}
=\varnothing,$$ then for any  $1\leq p\leq f_1-1$ there exists
$i_{(p)} = F_{1p}\setminus F_{1,p+1}$. Denote $j_{(p)} = c+f_1-p-1$.

Then for  $j=j_{(p)}$  we have either $w^{(j)}(i_{(p)})=j$, or
$w^{(j)}(i_{(p)})\in F_{1,p+1}$. In any case $(i_{(p)},
j_{(p)})\in S$, this proves  5). The statement  6) can be proved similarly. $\Box$\\
 {\bf Corollary 3.17}. \emph{All places   $(i,j)$, where $j\in J_a$ and $i>
 E_{[a]}$, are filled by zeros in the matrix $\Phi_\LC$ (resp. filled by symbols
"$\bullet$"{} in the diagram  $\DC$)}.

Denote by  $d_i$ and  $f_i$ the numbers of elements in  $D_i$ and
$F_i$
respectively. Let  $\nu$ be the greatest number such that  $D_{[\nu]}\subset J$. Denote $d_*=d_1+\ldots +d_\nu$.\\
 {\bf Proposition
3.18}. \emph{Let  $\xi$ satisfy the condition of case  2.  We claim that \\
 1) ~$\deg M_\xi(\la)= d_*$;\\
 2)~ $\deg M_\xi^\downarrow(\la)\leq d_*-1$, where $M_\xi^\downarrow(\la)$ is a minor that we get from
 $M_\xi(\la)$ moving one of the  rows on one line below.}\\
 {\bf Proof}.\\
 A) The minor  $M_\xi(\la)$,  as a polynomial in  $\la$, can be decomposed in the form
 $$M_\xi(\la)=\sum_{r=1}^{d} M_r \la^r, ~~\mbox{где}~~ d=\deg M_\xi(\la).\eqno(3.13) $$
 The coefficient $M_r$  is the sum of some minors of the matrix $\Phi_\LC$:
 $$ M_r = \sum_{R\subset I\cap J,~ |R|=r} M_R, ~~\mbox{где}~~
 M_R=M^{J\setminus R}_{I\setminus R}.\eqno(3.14) $$
A1) Show that  $d\leq d_*$. It is sufficient to prove that any minor
 $M_R$ with  $r=|R|>d_*$ equals to zero. Denote
$R_i=R\cap E_i$ and $r_i=|R_i|$. By assumption,
$$r=r_1+\ldots +r_\nu > d_1+\ldots+d_\nu=d_*.$$
Let  $a$ be a least number such that
$$\sum_{i=1}^a r_i
> \sum_{i=1}^a d_i.\eqno(3.15)$$ Note that  $a\leq \nu$.

We conclude the proof of  A1) in every of these cases separately:  i)~ $1\leq a < \nu$  or $a=\nu<l$ and ii)~ $a=\nu=l$.\\
 i)~ $1\leq a < \nu$  or $a=\nu<l$.
  By definitions  3.15,
$$|J_a| = \sum_{i=1}^a f_i + \sum_{j=1}^{b(a)} d_j,$$
$$|E_{[a]}| = \sum_{i=1}^a f_i + \sum_{j=1}^{a} d_j.$$

By lemma  3.16, we conclude that   $1\leq b(a)\leq a-1$. Hence,
$$ \sum_{i=1}^{b(a)} r_i\leq \sum_{i=1}^{b(a)} d_i.$$

Since  $J_a\subseteq E_{[b(a)]}$, then $J_a\cap R\subseteq
R_1\sqcup\ldots\sqcup R_{b(a)}$ and, therefore,
$$|J_a\setminus R| =|J_a|- |J_a\cap R|
\geq \sum_{i=1}^a f_i +   \sum_{i=1}^{b(a)} (d_i-r_i)\geq
\sum_{i=1}^a f_i.$$ Applying inequality  (3.15), we obtain
$$ |J_a\setminus R| \geq \sum_{i=1}^a f_i >  \sum_{i=1}^a f_i +   \sum_{i=1}^a
(d_i-r_i) = |E_{[a]}\setminus R|.\eqno(3.16)$$ By corollary 3.17,
the matrix  $\Phi_\LC$ has zeros in on all places  $(i,j)$, where $
j\in {J}_a$ and $i> E_{[a]}$. The minor  $M_R$ is zero, since all
its elements out of the rectangle  $\left(E_{[a]}\setminus R\right)
\times \left(J_a\setminus R\right)$, where $|E_{[a]}\setminus
R|>|J_a\setminus R|$,
are zero.\\
ii) ~ $a=\nu=l$.  In this case $E_{[a]}$ (resp. $J_a$) is defined in
 3.14(3) (resp. 3.14(5)). The statement of  A1) may be proved
 similarly to the case  i), changing  $f_l$ by  $f_{l1}=|F_{l1}|$.\\
A2). We conclude the proof of  A). It is necessary to show that
$M_{d*}\ne 0$.
  By formula  (3.14), any coefficient  $M_r$ is a sum of minors
$M_R$. The system of nonzero summands  $\{ M_R|~ |R|=r\}$ of this
sum is linear independent. Hence, the coefficient  $M_r\ne 0$ if and
only if there exists a minor $M_R\ne 0$ in the sum (3.14).

The coefficient   $M_{d_*}$ contains as a summand the minor
$$M_D = \prod _{i=1}^{l-1} M_{F_i}^{F_i'}\cdot M_{F_{l1}}^{F'_{l1}} \cdot M_{F_{l2}}^{F'_{l2}},\eqno(3.17)$$
that, by statements  5) and  6) of lemma  3.16, do not equal to
zero. Therefore,  $M_{d_*}\ne 0$ and $d=d_*$. This proves  A). \\
 B) Let as above  $I$ (resp. $J$) be the system of rows (resp. columns) of the minor  $M_\xi(\la)$. The system of columns for the minor
 $M_\xi^\downarrow(\la)$ the same as for  $M_\xi(\la)$, that is  $J$.
There exists  $g\in I$ such that  $g+1\notin I$ and the system of
rows $I^\downarrow$ for the minor  $M_\xi^\downarrow(\la)$ coincides
with  $(I\setminus \{g\})\cap\{g+1\}$.

By definition of the system  $I$, the number $g$ is the greatest in
some  $F_p$, where $\nu\leq p <l$. Then  $g+1\in D_{p+1}$. Similarly
to  $M_\xi(\la)$, the minor  $M_\xi^\downarrow(\la)$
 is presented in the  form
 $$M^\downarrow_\xi(\la)=\sum_{r=1}^{d(\downarrow)} M^\downarrow_r \la^r, ~~\mbox{where}~~
 d(\downarrow)=\deg M^\downarrow_\xi(\la). $$
 Each coefficient  $M^\downarrow_r$ is a sum of minors of the matrix  $\Phi_\LC$:
 $$ M^\downarrow_r = \sum_{R\subset I\cap J,~ |R|=r} M^\downarrow_R$$
 where  $  M^\downarrow_R$  is a minor of the matrix  $\Phi_\LC$ with system of rows
 $I^\downarrow\setminus R$ and system of columns  $J\setminus R$.
Let us show that any minor  $M^\downarrow_R$ equals to zero, if
$r\geq
 d_*$. This implies the last statement of proposition  3.18.

 Case  $r>d_*$ can be considered similar to   A1).
 Let  $r=d_*$. The case, when for some $1\leq a\leq\nu$ the inequality
   (3.15) is true,  also can be considered similar to  A1). We have to consider the case $d_i=r_i$ for any  $1\leq i\leq \nu$.

 Let  $I_p$ be the subsystem of rows that consists of  $i\in I$ such that  $i\leq g$. By definition of the system  $I$,
$$I_p= E_1\sqcup\ldots\sqcup E_\nu\sqcup F_{\nu+1}\sqcup\ldots\sqcup
F_p.$$

 As above we consider the subset  $J_p$ of columns  that is a union of all  $F'_i$,~ $1\leq i\leq p$, and all  $D_j$, ~ $1\leq j
\leq b(\nu)$ (see notation  3.15(4)). We obtain
$$|J_p\setminus R| \geq  \sum_{i=1}^p f_i + \sum_{i=1}^\nu d_i - \sum_{i=1}^\nu
r_i =\sum_{i=1}^\nu (f_i + d_i) +  \sum_{i=\nu+1}^p f_i -
\sum_{i=1}^\nu r_i = |I_p\setminus R|. \eqno(3.18)$$ By corollary
3.17,  all places  $(i,j)$, where $j\in J_p$ and $i>I_p$, are filled
by zeros in the matrix  $\Phi_\LC$. Respectively, in the minor $M_R$
all places  $(i,j)$, where $j\in J_p\setminus R$ and $i>I_p\setminus
R$ will be filled by zeros. Hence, in the minor  $M^\downarrow_R$
all places $(i,j)$, where $j\in J_p\setminus R$ and
$i>I_p\setminus(R\sqcup \{g\})$ will be filled by zeros. By (3.18),
$$|J_p\setminus R| > |I_p\setminus(R\sqcup \{g\})|.$$
Therefore, $M^\downarrow_R=0$.
$\Box$\\
{\bf Proposition  3.19}. \\\emph{
1) For any  $\xi\in S$ the minor  $M_{I,\LC}^J(\la)$ is extremal.\\
2) For any  $\xi\in S$ the element  $P_\xi$ of $K[LC^*]$  is
invariant with respect to the coadjoint representation of the group $L$.\\
3) Every  $P_\xi$ can be decomposed in the form  $P_\xi=y_\xi Q_\xi+
R_\xi$ where  $Q_\xi$ and $R_\xi$  belong to the subalgebra in
$S(\LC)$, generated by
$y_{ij}$,~ where $1\leq j< t$, and $ y_{it}$, where   $(i,t)\succ \xi$.} \\
{\bf Remark}. Since $P_\xi$ is invariant, then  $Q_\xi$ is also invariant.\\
{\bf Proof}. By the theorem 2.5, the statement  2) follows from the
statement 1).  In the case 1  the minor  $M_{\xi}(\lambda)$ has zero
degree. Since any root  $(w(j),j)$, where $j\in J$, belongs to $S$,
then $M_{\xi}(\lambda)= P_\xi\ne 0$. The lemma 3.6 implies that
 $M_\xi(\lambda)$ is extremal in the case 1.
By proposition 3.18 and lemma 3.6, $M_\xi(\lambda)$ is extremal in
the case 2.

In the case 1 one can prove the statement 3) decomposing the minor
 $P_\xi$ by its last $t$th column. To prove 3)  in the case 2, it is necessary to decompose by the last
 $t$th column all minors  $M_R$  in (3.14)  for   $r=d_*$ and $t\notin R$.\\
 {\bf Theorem 3.20}. \emph{The field of invariants of the coadjoint representation  of Lie algebra  $\LC$
 is isomorphic to the field of rational functions
 of  $P_\xi$, where
$\xi\in S$.}\\
{\bf Proof }. The theorem  1.1(2) implies that the field of
invariants of the coadjoint representation of $\LC$ is isomorphic to
the field of rational functions of $z_1,\ldots, z_s$. One can
conclude the proof of theorem, using the  induction by  $1\leq i\leq
s$ and applying the statements of theorem
1.1(1) and proposition  3.19(3).$\Box$ \\
{\bf Remark 3.21}. The polynomial  $z_i$ that was  constructed by
induction in the paper  ~\cite{P1} do not coincide with $P_{\xi_i}$
in general.
\\
{\bf Conclusion of the example  1}. The generators of the filed of
invariants of example 1 have the form
$$P_{\xi_1} = y_{41},\quad P_{\xi_2} = y_{62},\quad
P_{\xi_3} = y_{73},\quad P_{\xi_4} = y_{74}y_{41}+
y_{73}y_{31},\quad P_{\xi_5} = \left|\begin{array}{ccc}
y_{52}&y_{53}&y_{54}\\$$
y_{62}&y_{63}&y_{64}\\
0&y_{73}&y_{74}
\end{array}\right|.$$
For all  $\xi\in S$, but not for  $\xi_4$, the case 1 takes place,
and $P_\xi$ is a minor of  $\Phi_\LC$. In the case $\xi= \xi_4$, the
minor of the characteristic matrix
$$M_{\xi_4}(\la) =
 M^{1,2,3,4}_{2,3,4,7}(\la) = \left|\begin{array}{cccc}y_{21}& -\la&
0& 0\\
y_{31}&y_{32}& -\la&
0\\
y_{41}&y_{42}& y_{43}&-\la\\
0&0& y_{73}&y_{74}\end{array}\right|$$ is extremal
 and  $P_{\xi_4}$ is its
highest coefficient.

Samara State University\\
443011, Samara\\
ul. akad. Pavlova, 1,\\
Russia\\
\emph{ E-mail}: apanov@list.ru


\begin{thebibliography}{100}

\bibitem{K-Orb}
Kirillov A.A., \emph{ Lectures on the orbit method}, Graduate
Studies in Math. {\bf 64}(2004).
\bibitem{K-62}
Kirillov A.A., \emph{ Unitary reperetations of nilpotent Lie
groups},
 Uspekhi Mat. Nauk {\bf 17}(1962), 57-110.
 \bibitem{P1}
Panov A.N., \emph{ On index of certain nilpotent Lie algebras},
Contemp. math. and its applications {\bf 60}(2008), 123-131,
arXiv:0801.3025
\bibitem{P2}
Panov A.N., \emph{ Diagramm method in research on coadjoin orbits},
Vestnik SamGU, Natural Science Series (2008), no. 6(65), 139-151;
arXiv:0902.4584
\bibitem{P3}
Ignatev M.V., Panov A.N., \emph{ Coadjoint orbits for the group
$\mathrm{UT}(7,K)$}, Fundamental and applied math. {\bf 13}(2007),
no. 5, 127-159.
\bibitem{P4}
Panov A.N., \emph{Involutions in  $S_n$ and associated coadjoint
orbits},
 Zapiski POMI, 2007, {\bf 349}, 150-173, arXiv:0801.3022
\bibitem{Dix}
Dixmier J., \emph{ Algebres Enveloppantes}, Paris, 1974.
\bibitem{PW}
Parshall B., Wang J-P.,\emph{ Quantum linear groups}, Memoirs AMS
{\bf 89}(1991),  no. 439.
\bibitem{Ga}
Gavarini F., \emph{ Presentation by Borel subalgebras and Chevalley
generators for quantum enveloping algebras}, Proc. Edinburgh
Math.Soc. {\bf 49}(2006), 291-308.
\bibitem{Ja}
Jantzen J.G., \emph{ Lectures on Quantum Groups}, Graduate Studes in
Math. {\bf 6}(1995).

\end{thebibliography}
\end{document}